\documentclass[12pt,diff_eng]{article} 
\usepackage[english]{babel}
\usepackage{amssymb,amscd,amsmath}     
\begin{document}


\begin{center}
{\Large  
SETS OF INVARIANT MEASURES AND\\[5pt]
CESARO STABILITY}
\end{center}

\begin{center}{\large
\emph{S. G. Kryzhevich}
}\end{center}
 
\medskip  

\begin{center}
Saint-Petersburg State University, \\
Universitetskaya nab., 7-9,  Saint-Petersburg, Russia, 199034,\\
University of Nova Gorica, \\
Vipavska cesta, 13, Nova Gorica, Slovenia, SI-5000,\\
kryzhevicz@gmail.com. 
\end{center} 
 
\medskip

{\small 

\noindent\textbf{Abstract.} We take a space $X$ of dynamical 
systems that could be: homeomorphisms or continuous maps of a compact 
metric space $K$ or diffeomorphisms of a smooth manifold or actions of 
an amenable group. We demonstrate that a typical dynamical system of $X$ 
is a continuity point for the set of probability invariant measures 
considered as a function of a map, let $Y$ be the set of all such continuity points. 
As a corollary we prove that for dynamical systems of the residual set $Y$ average 
values of continuous functions calculated along trajectories do not drastically 
change if the system is perturbed. 
 
\medskip

}


\bigskip

\section{Introduction.} 

We consider the problem of structural stability of general dynamical systems that can be continuous maps/homeomorphisms of a metric compact space, diffeomorphisms of a smooth manifold or actions of a finitely generated group. We are interested in questions related to structural stability of such systems. For example, could we say that a generic dynamical system is structurally stable in a sense? If structural stability means topological conjugacy with close maps, then for $C^1$ diffeomorphisms of a smooth manifold, the answer is "no" \cite{smale}. However, individual stability of all periodic points is generic since Kupka -- Smale diffeomorphisms are generic \cite{kupkasmale}.  

An approach to study structural stability is given by the so-called Shadowing Theory \cite{shadowing} . Here both shadowing and inverse shadowing \cite{piljugincorless}, \cite{iranci} phenomena are important. Shadowing implies existence of exact trajectories in neighbourhoods of approximate ones and the inverse shadowing implies existence of trajectories of a fixed method in neighbourhoods of solutions. Observe that shadowing and inverse shadowing do not imply existence of a periodic solution in a neighbourhood of a periodic pseudotrajectory.  However, sometimes, shadowing implies robustness of periodic solutions. Sakai \cite{sakai} demonstrated that the $C^1$ - interior of the set of all diffeomorphisms with shadowing coincides with the set of all structurally stable diffeomorphisms. Osipov, Pilyugin and Tikhomirov \cite{opt} demonstrated that the so-called Lipschitz periodic shadowing property is equivalent to $\Omega$ -- stability, see also \cite{palmerpiljugintikhomirov}. Moreover, the corresponding set of dynamical systems coincides with the interior of the set of systems with periodic shadowing property and with the set of systems with orbital limit shadowing property. Pilyugin and Tikhomirov \cite{piltikh} demonstrated that Lipschitz shadowing is equivalent to structural stability. In all quoted cases periodic solutions are preserved. It was proved that shadowing is generic in the space of homeomorphisms of compact manifolds \cite{pipl}. However, it is not generic in the space of $C^1$ - diffeomorphisms \cite{bonattidiazturcat} see also \cite{yuanyorke} for $C^r$ - case.

F.\,Takens (see \cite{takens}, \cite{white}, \cite{mazur} and references therein) offered another approach, called tolerance stability. He proved that for a generic homeomorphism (or diffeomorphism) any trajectory of a perturbed map belongs to a neighbourhood of a trajectory of the non-perturbed one. We give precise formulation of corresponding theorems in Section 4. 

The aim of this paper is the following. It is proved that for a generic dynamical system  all points are orbit-wise weakly structurally stable. This means that for any trajectory of the initial system there is a trajectory of the perturbed system that lies in a neighbourhood of the initial trajectory.  Moreover, for an a priori selected observable, we can say that its average value along the given trajectory does not drastically change if the system is slightly perturbed. In a sense, the weak shadowing, we obtain, is close to statistical shadowing, introduced by M. Blank \cite{blank}, see also \cite{ow}. 

To prove this result methods of Ergodic Theory are applied. Also, the approach developed by author in \cite{vestnik} is used. First of all, we consider sets of Borel probability invariant measures and prove that a generic dynamical system is a "point of continuity"\ for such sets. 

Since this is true for measures, associated with trajectories, there must be a solution that engenders an invariant measure, close to initial one. This implies that there is a point such that its perturbed trajectory is close to the non-perturbed one "in average". This approach is very close to the Takens' theory. We discuss and compare results in Section 4.

\section{Stability of the set of invariant measures} 

Let $(K,\rho)$ be a compact metric set, consider the space of all continuous maps $X:=C^0(K\to K)$. Let $\mathrm{dist}_X$ be the standard metrics in $X$. Let ${\cal M}(K)$ be the set of all Borel probability measures on $K$ endowed with the topology of $*$ -- weak convergence. This topology can be engendered by the so called Kantorovich - Rubinstein distance. Let us recall one of equivalent definitions of this metrics
$$W_1(\mu,\nu)=\sup\left(\int_K f\, d\mu -\int_{K}\, f \, d\nu\right),$$
where the supremum is taken over the all $1$ - Lipschitz continuous functions $f:K\to {\mathbb R}$. Observe that the set ${\cal M}(K)$ becomes compact in that metrics. Later on, we use the same notion $W_1(\cdot,\cdot)$ for distance between a compact subset $P\subset {\cal M}(K)$ and a point $q\in {\cal M}(K)$. This distance turns to zero if and only if $q\in P$.

Recall the notion of Hausdorff metrics. Let given the space $K$ consider compact subsets $P,Q\subset {\cal M}(K,T)$. Then we can define
$$d_H(P,Q)=\max_{p\in P}W_1(p, Q)+\max_{q\in Q}W_1(q, P).$$
This distance turns to zero if and only if sets $P$ and $Q$ coincide.

For all metric spaces we deal with in this paper, we introduce following notions: $B_\varepsilon(x)$ stands for the $\varepsilon$ - ball centered at $x$ and $U_\varepsilon(A)$ stands for $\varepsilon$ - neighbourhood of a compact subset $A$.

Given a map $T\in X$ introduce the set ${\cal M}(K,T)$ of all Borel probability measures, invariant with respect to $T$. This set is also compact in the metrics $W_1$ and non-empty by Krylov -- Bogolyubov Theorem. 

Consider the map ${\cal M}: T \to {\cal M}(K,T)$ that is a map from $X$ to the set $Z$ of all compact subsets of ${\cal M}(K)$ endowed with the Hausdorff metrics $d_H$. 

\noindent\textbf{Definition 2.1.} We say that the map ${\cal M}$ is \emph{upper semicontinuous} at the point $T_0$ if for any $\varepsilon>0$ there exists a $\delta>0$ such that $\mathrm{dist}_X(T_0,T_1)<\delta$ implies 
${\cal M}(K,T_1) \subset U_\varepsilon({\cal M}(K,T_0))$.

Here we order subsets of $K$ by inclusion. As usually, "upper semicontinuity" means that that the value of the function cannot drastically increase in a neighbourhood of the point.

Recall that a set $Y\subset X$ is called \emph{residual} if it is a countable intersection of open dense subsets of $X$.

The main objective of this section is the following statement.

\noindent\textbf{Theorem 2.2.} \emph{For any compact metric space $K$ there exists a residual subset 
$Y\subset X$ such that any $T\in Y$ is a continuity point for the map $M$.}

We start with the following statement that is almost evident.

\noindent\textbf{Lemma 2.3.} \emph{The map ${\cal M}$ is upper semicontinuous.}

\noindent\textbf{Proof.} Fix a map $T\in X$. Consider a sequence of maps $T_k:X\to X$, converging to $T$. Let  $\mu_k$ be the sequence of Borel probability measures, invariant with respect to $T_k$. Suppose that the sequence $\mu_k$ $*$ -- wealky converges to $\mu$. Check that $\mu$ is invariant with respect to $T$ that is 
\begin{equation}\label{invt}
\int_K \varphi \, d\mu= \int_K \varphi \circ T\, d\mu, \qquad \forall \varphi\in C^0(K\to {\mathbb R}). 
\end{equation}

Fix a test function $\varphi\in C^0(K\to {\mathbb R})$: 
$$
\int\limits_K \varphi \, d\mu = \lim\limits_{k\to \infty} \int\limits_K \varphi \, d\mu_k=
\lim\limits_{k\to \infty} \int\limits_K \varphi \circ T_k \, d\mu_k=$$
\begin{equation}\label{e12}
\lim\limits_{k\to \infty} \int\limits_K \varphi \circ T \, d\mu_k+
\lim\limits_{k\to \infty} \int\limits_K (\varphi \circ T_k - \varphi \circ T) \, d\mu_k
\end{equation}
Since $T_k$ converges to $T$ uniformly, the second term in the right hand side of \eqref{e12} tends to zero. Since the sequence $\mu_k$ converges $*$ -- weakly to $\mu$, the first term in that right-hand side converges to $\int_K \varphi \circ T\, d\mu$. This proves equality \eqref{invt}. $\square$

To finish the proof of Theorem 2.2 it suffices to apply the following result first given in \cite{takens}, see also \cite[Corollary 11.1]{piljugin}. We provide this statement here reformulating it in our terms.

\noindent\textbf{Lemma 2.4 (Takens).} \emph{Let $X$ be a metric space, ${\cal M}:X \to Z$ be an upper semicontinuous $($or lower semicontinuous$)$ map of $X$ to $Z$ that is a set of all closed subset of a compact metric space. Then there is a residual set $Y\subset X$ such that any point of $Y$ is a point of continuity for the map $\cal M$.} 

Together with Lemma 2.3, the last statement finishes the proof of Theorem 2.2.

\noindent\textbf{Remark 2.5.} \emph{Instead of $C^0(K\to K)$ we may consider any subspace $X$ of that space with a topology that is not coarser than the topology induced by the metric of $C^0$ $($this means that convergence in the topology considered implies convergence in the standard topology of $C^0(K))$. For instance, we may consider the space of homeomorphisms of a compact metric space $K$. Or, we can take a $C^r$ smooth manifold $K$ $(r\ge 1)$ and consider a set of $C^r$ smooth diffeomorphisms of $K$ $($later on, whenever speaking of smooth maps with $r>0$, we assume that $K$ is a compact manifold$)$. A statement, similar to Theorem 2.2, is still valid in the considered case.}
  
We can also expand the obtained result to actions of finitely generated groups.  For any group $G$, we consider $C^r$ actions ($r\ge 0$) of that group that is the set of all homomorphisms $h:G\to \mathrm{Diff}^r(K)$ of $G$ to the set of $C^r$ smooth diffeomorphisms of $K$. If the group $G$ is amenable, the Krylov-Bogolyubov theorem is still applicable \cite{malpot} and the set of Borel invariant measures is non-empty.   

Fixing $\gamma=\{g_1,\ldots g_n\}$ that is the generating set of the group $G$, we introduce the distance $d_G$ between actions $h_1$ and $h_2$ by the formula
\begin{equation}\label{action}
d_\gamma(h_1,h_2)=\sum_{j=1}^n (d_r(h_1(e_j),h_2(e_j))+d_r(h_1(e_j^{-1}),h_2(e_j^{-1})))
\end{equation}
where $d_r$ in the right hand side of Eq. \eqref{action} stands for $C^r$ -- distance between maps.

\noindent\textbf{Corollary 2.6.} \emph{Let $r\in {\mathbb N} \bigcup \{0\}$. Take a compact metric space $K$ that is a $C^r$ -- smooth Riemannian manifold if $r>0$. Let $G$ be an amenable finitely generated group and $X$ be the set of all $C^r$ actions of $G$ with the metrics $d_G$ given by Eq.\, \eqref{action}. Consider the set $Z$ of all compact subsets in the set of all Borel probability measures on $K$ with the Hausdorff metrics engendered by Kantorovich - Rubinstein distance between measures. Given an action $h \in X$, consider ${\cal M}(h)$ that is the set of all Borel probability $h$ -- invariant measures. Then there exists a residual subset $Y\subset X$ such that any $h\in Y$ is a continuity point for the map $M$.}

One can get a result similar to Corollary 2.5, for actions of a finitely generated non-amenable groups (for instance, free groups). Then, we have to take the space $X$ of all actions that have a Borel probability invariant measure (e.g. actions of free groups that have a finite orbit).

\section{Cesaro structural stability.}

Let $r\in {\mathbb N}\bigcup \{0\}$. Consider a compact set $K$ that is a $C^r$ -- smooth manifold if $r>0$. Let $X$ be the set of all $C^r$ -- diffeomorphisms of $K$ (homeomorphisms if $r=0$). Fix the set $Y$ that exists by Theorem 2.2 or by Corollary 2.5. Given a point $p \in K$ and a map $T\in C^0(K\to K)$, we introduce a positive semiorbit $O^+_T(p):=\{T^k(p):k\ge 0\}$. If $T$ is  homeomorphism, we also consider the orbit of $p$: $O_T(p):=\{T^k(p):k\in {\mathbb Z}\}$. 

Next statement gives a "statistical" version of the so-called Zeeman's tolerance stability conjecture, see \cite{takens}.  

Given a function $\varphi\in C^0(K\to {\mathbb R})$ and a point $p$ we consider two values:
\begin{equation}\label{phipm}
\varphi_p^-=\liminf_{n\to \infty} \dfrac1{2n+1} \sum_{k=-n}^n \varphi(T^k(p)) \quad \mbox{and} 
\quad \varphi_p^+=\limsup_{n\to \infty} \dfrac1{2n+1} \sum_{k=-n}^n \varphi(T^k(p)).
\end{equation}

\noindent\textbf{Lemma 3.1.} \emph{Let $T\in X$ be a continuity point of the map $\cal M$, $\varphi\in C^0(K\to {\mathbb R})$. For any 
$\varepsilon,\sigma>0$ there exists $\delta>0$ such that the following statement is satisfied.  For any $p\in K$, and any $S\in X$ such that  
\begin{equation}\label{ledelta}
d_{C^0}(S,T)<\delta, \qquad d_{C^0}(S^{-1},T^{-1})<\delta
\end{equation}
there exists points $q_\pm$ such that
\begin{equation}\label{statstab0}
\liminf_{n\to \infty} \dfrac{\# \{k=-n,\ldots,n : S^k(q_\pm)\in U_\sigma (\overline{O_T(p)})\}}{2n+1} \ge 1-\varepsilon;
\end{equation}
and there exists  limit
\begin{equation}\label{statstab1}
\lim_{n\to \infty} \dfrac1{2n+1} \sum_{k=-n}^n\varphi(S^k(q_-))\ge \varphi_p^- -\varepsilon, \qquad
\lim_{n\to \infty} \dfrac1{2n+1} \sum_{k=-n}^n\varphi(S^k(q_+))\le \varphi_p^+ +\varepsilon.
\end{equation}}

\noindent\textbf{Proof.} \textbf{Step 1.} Fix a map $T$,  positive values $\varepsilon<1$ and $\sigma$ and the function $\varphi$. Without loss of generality, we may assume that $\sup_K |\varphi(x)|\le 1$. 

Fix a value $\alpha \in (0,1/2)$, now this value is not important, we need it later. Fix a $C^\infty$ smooth function $\eta: {\mathbb R}\to {\mathbb R}$ such that
$$\begin{array}{rl}
\eta(t)=1 \qquad & \mbox{for all} \qquad t \le 1-\alpha;\\
\eta (t)=0 \qquad & \mbox{for all} \qquad t \ge 1;\\
\eta'(t) \in (-2/\alpha,0) & \mbox{for all} \qquad t \in (1-\alpha,1).
\end{array}$$
Consider a function $\psi(x)=\eta(\rho(x,{\overline{O_T(p)}})/\sigma)$. Observe that this function is Lipschitz continuous with the constant $L_0:=2/(\alpha\sigma)$.

There exists an increasing sequence $\{n_m\in {\mathbb N}\}$ and a Borel $T$ -- invariant Borel probability measure $\mu_T$ such that for any continuous function $f:K\to {\mathbb R}$ we have 
\begin{equation}\label{mut}
\lim_{n\to \infty} \dfrac1{n_m} \sum_{k=0}^{n_m-1} f(T^k(p))=\int_K f\, d\, \mu_T
\end{equation}
(see \cite[Theorem 4.4.1]{kaha}). Evidently, $0\le \psi(x)\le 1$ for all $x$ and
$$\int_K \psi\, d\mu_T =1.$$

Take $\delta>0$ so small that for any $S$, satisfying Eq. \eqref{ledelta}, there is an $S$ -- invariant measure $\mu_S$ such that 
\begin{equation}\label{psis}
\int_K \psi\, d\mu_S \ge 1-\dfrac{\varepsilon^2}8.
\end{equation}
This can be done by Corollary 2.5. 

By two-side Birkhoff Ergodic Theorem \cite[Theorem 4.1.3]{kaha} there exists a function ${\hat \psi} \in {\mathbb L}^1(\mu_S)$ that equals to 
$$\lim_{n\to \infty}\dfrac1{2n+1}\sum_{k=-n}^n \psi \circ S^k$$
$\mu_S$ almost everywhere. Observe that 
$$\int_K \psi \, d\mu_S= \int_K {\hat \psi}\, d\mu_S.$$

\noindent\textbf{Proposition 3.2.} \emph{For any measurable set $A$ such that $\mu_S(A)\ge \varepsilon/4$ there exists a point $q\in A$ such that ${\hat \psi}(q)> 1-\varepsilon$.}

\noindent\textbf{Proof.} Recall that $\psi\le 1$ almost everywhere w.r.t. $\mu_S$. If the statement of the lemma is wrong, we have 
$$\int_K {\hat\psi}\, d\mu_S = \int_A {\hat\psi} \, d\mu + \int_{K\setminus A} {\hat\psi} \, d\mu\le (1-\varepsilon)\cdot \varepsilon /4+1\cdot (1-\varepsilon/4)=1-\varepsilon^2/4$$
that contradicts to  Eq. \eqref{psis}. $\square$

Observe that the function $\psi$ vanishes out of the set $U_\sigma(\overline{O_T(p)})$ which demonstrates that any set $A$: $\mu_S(A)\ge \varepsilon$ contains a point $q$ such that the formula similar to \eqref{statstab0} is satisfied.

\textbf{Step 2.} Suppose that $\delta$ is so small that if Eq.\, \eqref{ledelta} is satisfied then $\mu_S$, satisfying requirements of Step 1, could be taken so that 
\begin{equation}\label{ST}
\left| \int_K \varphi\, d\mu_S - \int_K \varphi\, d\mu_T\right|\le \varepsilon/4.
\end{equation}
 
Now let us prove existence of the point $q_+$ satisfying conditions of the theorem. Existence of the point $q_-$ can be proved similarly. Apply the two-side Birkhoff Ergodic Theorem to the function $\varphi$. There exists a function 
${\hat\varphi}\in {\mathbb L}^1(\mu_S)$ that equals 
$$\lim_{n\to \infty} \dfrac1{2n+1} \sum_{k=-n}^n\varphi \circ S^k$$
almost everywhere. Also we have 
$$\int_K \hat\varphi\, d\mu_S=\int_K \varphi\, d\mu_S$$
and $|\hat\varphi(x)|\le 1$ almost everywhere. It follows from equations \eqref{phipm}, \eqref{mut} and \eqref{ST} that
$$\int_K \hat\varphi \, d\mu_S \ge \varphi_p^--\varepsilon/4.$$ 
Let $A=\{x\in K: \hat\varphi(x)\ge  \varphi_p^--\varepsilon\}$. Then $\mu_S(A)\ge \varepsilon/4$. Otherwise, since ${\hat\varphi}\le 1$ a.e., $\varphi_p^-\ge -1$ and $\varepsilon<1$ we have
$$\int_K {\hat\varphi}\, d\mu_S = \int_A {\hat\varphi} \, d\mu + \int_{K\setminus A} {\hat\varphi} \, d\mu\le \dfrac{\varepsilon}4+\left(1-\dfrac{\varepsilon}4\right)\left(\varphi_p^--\varepsilon\right)<
\varphi_p^--\dfrac{\varepsilon}4$$
that gives a contradiction. Hence, by Lemma 3.2 there exists a point $q_+$ that satisfies requirements of the theorem. $\square$ 

Let us list some versions of Theorem 3.1. The next one concerns actions of non-invertible continuous maps, so this is a one-side version of Theorem 3.1. We formulate it for 
$C^0$ -- maps of a compact set to itself but it can also be formulated for $C^r$ -- maps.

\noindent\textbf{Theorem 3.3.} \emph{Let $X=C^0(K\to K)$, $Y$ be the set of continuity points of the map ${\cal M}$ (that is residual by Theorem 2.2). Take $T\in Y$, $\varphi\in C^0(K\to {\mathbb R})$. For any 
$\varepsilon, \sigma>0$ there exists $\delta>0$ such that the following statement is satisfied.  For any $p\in K$, and any $S\in X$ satisfying Eq.\, \eqref{ledelta}
there exist two points $q_-$ and $q_+$ such that
$$
\liminf_{n\to \infty} \dfrac{\# \{k=0\ldots,{n-1} : S^k(q_\pm)\in U_\sigma (\overline{O_T(p)})\}}{n} \ge 1-\varepsilon;
$$
and there exist limits
$$
\lim_{n\to \infty} \dfrac1{n} \sum_{k=0}^{n-1}\varphi(S^k(q_+))\ge \varphi_p^- -\varepsilon, 
\qquad \lim_{n\to \infty} \dfrac1{n} \sum_{k=0}^{n-1} \varphi(S^k(q_-))\le \varphi_p^+ +\varepsilon.
$$}

The proof is similar to proof of Theorem 3.1 with the only difference: we use one-side Birkhof Ergodic Theorem \cite[Theorem 4.1.2]{kaha}.

Basing on Corollary 2.6, we could prove a similar result for actions of finitely generated amenable groups. Such actions could be treated as actions of finitely generated free groups where the existence of Borel probability invariant measures is guaranteed. The Birkhoff's theorem for actions of free groups was provided by Stein and Nevo \cite{nest}.

In the rest of this section we always deal with homeomorphisms of compact metric spaces having in mind all possible generalisations.

For Lipschitz continuous functions $\varphi$, we have the following "uniform"\ version of Theorem 3.1:

\noindent\textbf{Theorem 3.4.} \emph{For all $T\in Y$ and all 
$\varepsilon, \sigma, L>0$ there exists $\delta>0$ such that the following statement is satisfied.  For any $L$ -- Lipschitz continuous function $\varphi: K\to {\mathbb R}$, any point $p\in K$, and any $S\in X$ satisfying Eq.\, \eqref{ledelta} there exist two points $q_{\varphi}^-$ and $q_{\varphi}^+$ such that
$$
\liminf_{n\to \infty} \dfrac{\# \{k=-n,\ldots,n : S^k(q_\varphi^\pm) \in U_\sigma (\overline{O_T(p)})\}}{2n+1} \ge 1-\varepsilon;
$$
and there exist limits
$$ \lim_{n\to \infty} \dfrac1{2n+1} \sum_{k=-n}^n\varphi(S^k(q_{\varphi}^+))\ge \varphi_p^- -\varepsilon, 
\qquad \lim_{n\to \infty} \dfrac1{2n+1} \sum_{k=-n}^n \varphi(S^k(q_{\varphi}^-))\le \varphi_p^+ +\varepsilon.
$$}

We proof repeats one of Theorem 3.1 due to definition of Kantorovich metric $W_1$.

Finally, we give a result on visiting a domain by iterations of a trajectory. Let $V\subsetneq K$ be an open nonempty set,  $\partial V$ be the boundary of $V$. 

Given $\beta>0$, we define two open sets: 
$$V_+^{\beta}=U_{\beta}(\overline V), \qquad  
V_-^\beta=V\setminus \overline{U_{\beta}(K\setminus V)}.$$
Introduce two Lipschitz continuous functions $\chi^\beta_+$ and $\chi^\beta_-$ such that
$$ \begin{array}{c}
\chi^\beta_+|_{V}=1, \qquad \chi^\beta_-|_{V_-^\beta}=1;\\
\chi^\beta_+|_{K\setminus V_+^{\beta}}=0, \qquad \chi^\alpha_-|_{K\setminus V}=0;\\
\chi^\beta_+(x)= \eta\left(\dfrac{\rho(x, \partial V)}{\alpha}\right), \qquad
\mbox{for all } \quad x\in V_+^{\beta} \setminus V;\\
\chi^\beta_-(x)= \eta\left(\dfrac{\rho(x, \partial V)}{\alpha}\right), \qquad
\mbox{for all } \quad x\in V \setminus V_-^{\beta};
\end{array}$$
where $\eta$ is the function defined in the proof of Theorem 3.1. Finally, we introduce the function 
$\chi(x)$ that is the characteristic function of the set $V$. Evidently, for any $\beta>0$ we have
$$\chi^\beta_+\ge \chi \ge \chi^\beta_-$$
everywhere in $K$.

Given a point $p\in K$ and an open set $V\subset K$, similarly to Eq.\, \eqref{phipm} we introduce values
$$
\chi_p^-=\liminf_{n\to \infty} \dfrac1{2n+1} \sum_{k=-n}^n \chi(T^k(p)) \quad \mbox{and} 
\quad \chi_p^+=\limsup_{n\to \infty} \dfrac1{2n+1} \sum_{k=-n}^n \chi(T^k(p)).
$$
Applying Theorem 3.1 to functions $\chi^\beta_+$ and $\chi^\beta_-$, we get the following statement.

\noindent\textbf{Theorem 3.5.} \emph{For all $T\in Y$, any open set $V$ and all 
$\varepsilon, \sigma, \beta>0$ there exists $\delta>0$ such that the following statement is satisfied.  For any point $p\in K$, and any $S\in X$ satisfying Eq.\, \eqref{ledelta} there exist two points 
$q_{V}^-$ and $q_{V}^+$ such that
$$
\liminf_{n\to \infty} \dfrac{\# \{k=-n,\ldots,n : S^k(q_V^\pm) \in U_\sigma (\overline{O_T(p)})\}}{2n+1} \ge 1-\varepsilon;
$$
and there exist limits
$$ \lim_{n\to \infty} \dfrac{\# \{k=-n,\ldots,n : S^k(q_{V}^+)\in V\}}{2n+1}\ge \chi_p^- -\varepsilon, $$
$$ \lim_{n\to \infty} \dfrac{\# \{k=-n,\ldots,n : S^k(q_{V}^-)\in V\}}{2n+1}\le \chi_p^- +\varepsilon. $$} 

\section{Takens' theory.} 

In this section we compare our results with closely related statements of Takens' Theory \cite{takens}, \cite{white}, \cite{mazur}, \cite{piljugin} (see also references therein) and compare them with the result of our paper. 

\noindent\textbf{Definition 4.1.}  \emph{Let $(K,\rho)$ be a compact metric space, $H(K)$ be the set of all homeomorphisms of $K$. Fix $\varepsilon > 0$. We say that two homeomorphisms $S$ and $T$ are orbitally $\varepsilon$ -- equivalent if the following conditions are satisfied:
\begin{itemize}
\item[1.] For any trajectory $O_S(p)$ of the homeomorphism $S$ there exists a trajectory $O_T(q)$ of the homeomorphism $T$ such that 
\item[1.1.]  $O_S(p)\subset U_{\varepsilon} (\overline{O_T(q)})$,
\item[1.2.]  $O_T(q)\subset U_{\varepsilon} (\overline{O_S(p)})$;
\item[2.] for any trajectory $O_T(q)$ of the homeomorphism $T$ there exists a trajectory $O_S(p)$ of the homeomorphism $S$ such that 
\item[2.1.]  $O_T(q)\subset U_{\varepsilon} (\overline{O_S(p)})$,
\item[2.2.] $O_S(p)\subset U_{\varepsilon} (\overline{O_T(q)})$.
\end{itemize}}

\noindent\textbf{Definition 4.2.}  \emph{Let $D\subset H(M)$ be a subset with the inherited topology. 
We say that a homeomorphism $T\in D$ is tolerance $D$ -- stable if for any $\varepsilon>0$ there exists a neighbourhood $V$ of $T$ in $D$ such that any diffeomorphism $S\in V$ is $\varepsilon$ -- equivalent to $T$.}

F.\, Takens \cite{takens} formulated the following conjecture.

\noindent\textbf{Conjecture 4.3.} (Zeeman's tolerance stability conjecture). \emph{For any $D \subset H(M)$ there exists a residual $($in the topology of $D)$ subset $D_0\subset D$ such that any homeomorphism $T\in D_0$ is tolerance $D$ -- stable.}

In this form the statement is wrong, a counterexample was provided by W. White in \cite{white} (see also \cite{piljuginc0}) proved that generally, speaking, tolerance stability conjecture is wrong. However, Takens proved a weaker version of the conjecture. 

He considered two weakened variants of the property of tolerance $D$-stability. As above, fix $\varepsilon > 0$. We say that two homeomorphisms $f$ and $g$ are minimally $\varepsilon$ -- equivalent (maximally $\varepsilon$ -- equivalent) if conditions 1.1 and 2.1 of Definition 4.1 (respectively, conditions 1.2 and 2.2) are omitted in the above definition of orbital $\varepsilon$ -- equivalence.

Denote by $D_{max}$ the subset of $D$ which consists of homeomorphisms $T$ having the following property: For any $\varepsilon> 0$ there exists a neighbourhood $V$ of $T$ in $D$ such that any homeomorphism $S\in V$ is maximally $\varepsilon$ -- equivalent to $f$ (recall that we consider the set $D$ with a topology that is not coarser than the standard $C^0$ topology).

Similarly, we denote by $D_{min}$ the subset of $D$ which consists of homeomorphisms
$f$ having the following property: For any $\varepsilon > 0$ there exists a neighbourhood $V$ of $T$
in $D$ such that any homeomorphism $S \in  W$ is minimally $\varepsilon$ - equivalent to $T$.
Takens proved the following statement.

\noindent\textbf{Theorem 4.4 (Takens)} \emph{Any of the sets $D_{max}$ and $D_{min}$ is residual in $D$.}

A. Mazur \cite{mazur} demonstrated that for the complete space of homeomorphisms of a compact manifold Conjecture 4.3 is valid in a stronger form. Before formulating the statement we give one more definition. 

\noindent\textbf{Definition 4.5.} A homeomorphism $T\in H(M)$ is strongly tolerance stable if for every $\varepsilon > 0$ there exists $\delta > 0$ such that for every $S\in B_\delta(T)$ each $S$ -- orbit is $\varepsilon$ -- traced (see Definition 3.6) by some $T$ -- orbit and each $S$ --orbit is $\varepsilon$ -- traced by some $T$ -- orbit.

\noindent\textbf{Theorem 4.6 (Mazur).} \emph{Let $K$ be a compact smooth manifold with the metric d induced by the Riemanian structure. A generic $f \in H(K)$ has the strong tolerance stability property.}

For periodic points we can say a bit more.

\noindent\textbf{Definition 4.7.} \emph{Let $\varepsilon>0$, $S$ be a homeomorphism of a metric compact set $(K,\rho)$. We say that a sequence $\{x_i: i\in {\mathbb N}\}$ is $\varepsilon$ -- traced by a trajectory $O_S(p)$ if for any $i$ we have
$\rho (x_i,S^i(p))<\varepsilon$.}

\noindent\textbf{Theorem 4.8.} \emph{Let $X$ be a subset of $H(K)$ -- the space of all homeomorphisms of $K$ with the topology not coarser than one, inherited from $H(K)$. Suppose that any homeomorphism $T\subset X$ has a fixed point $($or a point with a period less or equal than a given number $n_0)$. Then there exists a residual set $Y$ such that for any $T\in Y$, any $\varepsilon>0$ and any $n\in N$ $(n\ge n_0)$ there exists a neighbourhood $V_n$ of $T$ in $X$ such that for any $S\in V_n$ and any $m_p$ -- periodic point $p$ of $T$ $(m\le n)$ there is a periodic point $q\in K$ with a period $m_q\le n$ such that the orbit $O_T(p)$ is $\varepsilon$ -- traced by $O_S(q)$.}

To prove this statement, it suffices to consider a map of the set $X$ to the set of all point with period $\le n$. This map is upper semicontinuous so we could apply Takens' result (Lemma 2.4) similarly to the proof of Theorem 2.2.

Observe that, for instance typical diffeomorphism of a smooth manifold has periodic points, as it was proved by Bonatti and Crovisier \cite{boncrov}. Moreover, they demonstrated that for a typical diffeomorphism the set of periodic points is dense in the set of all chain recurrent points.

Let us compare results of Theorems 4.4 and 4.6 with Theorems 3.1 -- 3.5. Theorem 4.4 works for a very wide class of dynamical systems and, unlike our formula \eqref{statstab0}, it claims that the whole set of one trajectory (not the majority of points) belongs to a neighbourhood of another one. However, only closeness of sets is provided and no "statistical"\ or "Cesaro average"\ properties like \eqref{statstab1} can be provided or derived strictly from the statement of Takens' theorem. 

For homeomorphisms of manifolds, Mazur's statement (one of Theorem 4.6) is much stronger than ones of our results. Pointwise $\varepsilon$ -- tracing is much more than all statistical properties, we are studying. However, the proof, provided in \cite{mazur} and, more generally, ideas of the paper \cite{pipl}, crucially depend on the manifold structure of the considered set. Moreover, the proof uses the fact that the complete set of homeomorphisms is considered, not a subset. For instance, it cannot be easily transformed to a statement on $C^r$ -- diffeomorphisms ($r\ge 1$).

We believe that our result stands in the middle between Takens' and Mazur's ones and thus completes the Takens' tolerance stability theory. 

\bigskip

\noindent\textbf{Acknowledgements}. This work was supported by Russian Foundation for Basic Researches, grant 15-01-03797-a. Author is grateful to Prof. Sergei Pilyugin for attracting attention to  tolerance stability approach.

\end{document}